\documentclass[12pt]{article}

\usepackage[centertags]{amsmath}
\usepackage{amsfonts}
\usepackage{amsthm}
\usepackage{newlfont}
\usepackage{amscd}
\usepackage{amsgen}
\usepackage{amssymb}
\usepackage{amssymb,amsmath}
\usepackage{amscd}
\usepackage{enumerate}
\usepackage{bbm}

\newlength{\defbaselineskip} 
\theoremstyle{plain}
\newtheorem{thm}{Theorem}[section]
\newtheorem{cor}[thm]{Corollary}
\newtheorem{con}[thm]{Conjecture}
\newtheorem{df}[thm]{Definition}
\newtheorem{lema}[thm]{Lemma}
\newtheorem{obs}[thm]{Proposition}
\newtheorem{exm}[thm]{Example}

\newtheorem{rem}[thm]{Remark}

\theoremstyle{definition} 
\theoremstyle{definition} \newtheorem{ex}{Example}[section] %

 \numberwithin{equation}{section}
\def\p{\mathbb{P}}
\def\r{\mathbb{R}}
\def\z{\mathbb{Z}}

\def\c{\mathbb{C}}
\def\q{\mathbb{Q}}
\def\o{\mathcal{O}}

 \DeclareMathOperator{\Ext}{Ext}
 
\DeclareMathOperator{\Pic}{Pic}
 
\DeclareMathOperator{\Spec}{Spec}
\DeclareMathOperator{\Hom}{Hom}
\def\p{\mathbb{P}}

\def\ob{\begin{obs}}
\def\kob{\end{obs}}
\def\dow{\begin{proof}}
\def\kdow{\end{proof}}

\def\tw{\begin{thm}}
\def\ktw{\end{thm}}
\def\hip{\begin{con}}
\def\khip{\end{con}}
\def\lem{\begin{lema}}
\def\klem{\end{lema}}
\def\ex{\begin{exm}}
\def\prog{\begin{pr}}
\def\kprog{\end{pr}}
\def\wn{\begin{cor}}
\def\kwn{\end{cor}}
\def\uwa{\begin{rem}}
\def\kuwa{\end{rem}}
\def\kex{\end{exm}}
\def\dfi{\begin{df}}
\def\kdfi{\end{df}}

\begin{document}

\title{On the full, strongly exceptional collections on toric varieties with Picard number three}


\author{Micha\l \hskip 5pt Laso\'{n}
\and
Mateusz Micha\l ek}


\maketitle
\begin{abstract}
We investigate full strongly exceptional collections on smooth, complete toric varieties. We obtain explicit results for a large family of varieties with Picard number three, containing many of the families already known. We also describe the relations between the collections and the split of the push forward of the trivial line bundle by the toric Frobenius morphism.
\end{abstract}

\tableofcontents

\section{Introduction}\label{sIntro}
Let $X$ be a smooth variety over an algebraically
closed field $\mathbb K$ of characteristic zero and let $D^b(X)$
be the derived category of bounded complexes of coherent sheaves
of $\mathcal O_X$-modules. This category is an important algebraic
invariant of $X$. In order to understand the derived category
$D^b(X)$ one is interested in knowing a  strongly exceptional
collection of objects that generate $D^b(X)$, see also \cite{bond}.

For a smooth, complete toric variety $X$ there is a well known
construction due to Bondal which gives a full collection of line
bundles in $D^b(X)$. In some cases Bondal's collection of line
bundles is a strongly exceptional collection (see also \cite{ober}), but it is not true in general. Often one can find a subset of this collection and order it in such a way that it becomes strongly exceptional and remains full. This approach was well described in \cite{lcmr} for a class of toric varieties with Picard number three.

One of the first conjectures concerning this topic was made by A. King \cite{king}: \hip[King's] For any smooth, complete
toric variety $X$ there exists a full, strongly exceptional
collection of line bundles.\khip

Originally this conjecture
was made in terms of existence of titling bundles whose direct
summands are line bundles, but it is easy to see that they are
equivalent, see \cite{lcmr2}.
It was disproved by Hille and Perling, in \cite{hipe}.
They gave an example of a smooth, complete toric surface which
does not have a full, strongly exceptional collection of line
bundles. The conjecture was reformulated by Mir\'o-Roig and Costa (stated also in \cite{bohu}):

\hip For any smooth, complete Fano toric variety there exists a
full, strongly exceptional collection of line bundles. \khip

This conjecture is still open and is supported by many numerical
evidence. It has an affirmative answer when the Picard number of $X$
is less then or equal to two \cite{lcmr2} or the dimension of $X$ is at most two \cite{bohu}. Recently it was also prooved for dimension three \cite{bond}, \cite{BeTi}. Even when the Picard number is equal to $3$ the question remains open.

The goal of this paper is to investigate when it is possible to
find a full, strongly exceptional collection and whether line bundles
that come from Bondal's construction contain such a collection. We restrict our attention to
smooth, complete toric varieties with Picard number three. There are some families among these varieties for
which the conjecture is true \cite{dlm}, \cite{lcmr}. We prove it
in section \ref{yes} for a greater family of varieties containing
both families already known. In section \ref{no}  we also show that in general it is not possible for a
smooth, complete toric variety with Picard number three to find a
full, strongly exceptional collection among line bundles that come
from Bondal's construction, even in the Fano case.

To determine the image of Bondals construction we look at the image of the real torus in the Picard group of a toric variety. We also compare this with the result of Thomsen's algorithm \cite{thom} that gives a decomposition of the push forward of a line bundle by a toric Frobenius morphism. This leads to some unexpected results like Corollary \ref{dz}.

To prove that a given collection of line bundles is strongly exceptional we develop new, efficient methods of counting homologies of simplicial complexes given by primitive collections, that is minimal subsets of points that do not form a simplex. To do this we use the results of \cite{mroz}. In particular this enables us to determine all acyclic simplicial complexes arising from complete toric varieties with Picard number three.

\section*{Acknowledgements}
We would like to thank very much Rosa Maria Mir\'o-Roig and Laura Costa for introducing us to this subject. We are also grateful for many useful and interesting talks.

\section{Preliminaries}\label{sPrelim}

\subsection{Full, strongly exceptional collections}\label{ssFullStrExColl}
For an algebraic variety $X$ let $D^b(X)$ be the derived category of coherent sheaves on $X$. For an introduction to derived categories the reader is advised to look in \cite{caladarau} and \cite{ganfeldmanin}. The structure and properties of the derived category of an arbitrary variety $X$ can be very complicated and they are an object of many studies. One of the approaches to understand the derived category uses the notion of exceptional objects. Let us introduce the following definitions (see also \cite{goru}):
\begin{df} \hskip 0pt
\begin{enumerate}
\item{} A coherent sheaf $F$ on $X$ is {\bf exceptional} if
$\text{Hom}(F,F)= \mathbb K$ and \text{Ext}$\,^i_{\o_X}(F,F)
\,=\,0$ for $i \geq 1$. \item{} An ordered collection
$(F_0,\,F_1,\,\dots,\,F_{m})$ of coherent sheaves on $X$ is an
{\bf exceptional collection} if each sheaf $F_i$ is exceptional
and\newline $\text{Ext}^{\hskip 2 pt i}_{\o_X}(F_k,F_j)=0$ for $j<k$ and $i\geq 0$.
\item{} An exceptional collection $(F_0,\,F_1,\,\dots,\,F_{m})$
of coherent sheaves on $X$ is a {\bf strongly exceptional
collection} if $\text{Ext}^{\hskip 2 pt i}_{\o_X}(F_j,F_k)=0$ for $j \leq k$ and
$i \geq 1$. \item{} A (strongly) exceptional collection
$(F_0,\,F_1,\,\dots,\,F_{m})$ of coherent sheaves on $X$ is a
{\bf full, (strongly) exceptional collection } if it generates the
bounded derived category $D^b(X)$ of $X$ i.e. the smallest
triangulated category containing $\{F_0,\,F_1,\,\dots,\,F_n\}$ is
equivalent to $D^b(X)$.
\end{enumerate}
\end{df}
For an exceptional collection $(F_0,\dots,F_m)$ one may define an object $F=\oplus_{i=0}^m F_i$ and an algebra $A=\Hom (F,F)$. Such an object gives us a functor $G_F$ from $D^b(X)$ to the derived category $D^b(A-mod)$ of right finite-dimensional modules over the algebra $A$. Bondal proved in \cite{bond}, that if $X$ is smooth and $(F_i)$ is a full, strongly exceptional collection, then the functor $G_F$ gives an equivalence of these categories. For further reading only the definition of the strongly exceptional collection is necessary.

\subsection{Toric varieties}\label{ssToricVar}
A normal algebraic variety is called toric in it contains a dense torus $(\c^*)^n$ whose action on itself extends to the action on the whole variety. For a good introduction to toric varieties the reader is advised to look in \cite{cox} or \cite{fult}. Varieties of this type form a sufficiently large class among normal varieties to test many hypothesis in algebraic geometry. Many invariants of a toric variety can be effectively computed using combinatorial description. Let us recall it.

Given an $n$ dimensional torus $T$ we may consider one parameter subgroups of $T$, that is morphisms $\c^*\rightarrow T$ and characters of $T$, that is morphisms $T\rightarrow\c^*$. One parameter subgroups form a lattice $N$ and characters form a lattice $M$. These lattices are dual to each other and isomorphic to $\z^n$.

A toric variety $X$ is constructed from a fan $\Sigma$, that is a system of cones $\sigma_i\subset N$. This is done by gluing together affine schemes $\Spec(\c[\sigma_i^*])$, where $\sigma_i^*\subset M$ is a cone dual to $\sigma_i$.
One dimensional cones in $\Sigma$ are called rays. The generators of these semigroups are called ray generators.

Many properties of the variety $X$ can be described using the fan $\Sigma$. For example $X$ is smooth if and only if for every cone $\sigma_i$ the set of its ray gene\-rators can be extended to the basis of $N$. Moreover to each ray generator $v$ we may associate a unique $T$ invariant Weil divisor denoted by $D_v$. There is a well known exact sequence:
\begin{align}\label{PPic}
0\rightarrow M\rightarrow Div_T\rightarrow\ Cl(X) \rightarrow 0,
\end{align}
where $Div_T$ is the group of $T$ invariant Weil divisors and $Cl(X)$ is the class group. The map $M\rightarrow Div_T$ is given by:
$$m\rightarrow \sum m(v_i)D_{v_i},$$
where the sum is taken over all ray generators $v_i$.

Smooth, complete
toric varieties with Picard number three have been classified by
Betyrev in \cite{baty} according to their primitive relations. Let
$\Sigma$ be a fan in $N=\z^n$.

\dfi We say that a subset $P\subset R$ is a primitive collection
if it is a minimal subset of $R$ which
does not span a cone in $\Sigma$. \kdfi

In other words a primitive collection is a subset of ray
generators, such that all together they do not span a cone in
$\Sigma$ but if we remove any generator, then the rest spans a
cone that belongs to $\Sigma$. To each primitive collection
$P=\{x_1,\dots, x_k\}$ we associate a primitive relation. Let
$w=\sum_{i=1}^k x_i$. Let $\sigma\in\Sigma$ be the cone of the
smallest dimension that contains $w$ and let $y_1,\dots,y_s$ be
the ray generators of this cone. The toric variety of $\Sigma$ was
assumed to be smooth, so there are unique positive integers
$n_1,\dots, n_s$ such that
$$w=\sum_{i=1}^s n_iy_i.$$
\dfi For each primitive collection $P=\{x_1,\dots,x_k\}$ let $n_i$
and $y_i$ be as described above. The linear relation:
$$x_1+\dots+x_k-n_1y_1-\dots-n_sy_s=0$$
is called the primitive relation (associated to $P$). \kdfi Using
the results of \cite{grun} and \cite{odap} Batyrev proved in
\cite{baty} that for any smooth, complete $n$ dimensional fan with $n+3$
generators its set of ray generators can be partitioned into $l$ non-empty
sets $X_0,\dots,X_{l-1}$ in such a way that the primitive
collections are exactly sums of $p+1$ consecutive sets $X_i$ (we
use a circular numeration, that is we assume that $i\in\z/l\z$), where
$l=2p+3$. Moreover $l$ is equal to $3$ or $5$. The number $l$ is
of course the number of primitive collections. In the case $l=3$
the fan $\Sigma$ is a splitting fan (that is any two primitive
collections are disjoint). These varieties are well characterized,
and we know much about full, strongly exceptional collections of
line bundles on them. The case of five primitive
collections is much more complicated and is our object of study. For
$l=5$ we have the following result of Batyrev \cite{baty}, Theorem 6.6:

\tw\label{charBatyrev} Let $Y_i=X_i\cup X_{i+1}$, where $i\in\z/5\z$,
$$X_0=\{v_1,\dots,v_{p_0}\},\quad X_1=\{y_1,\dots,y_{p_1}\},\quad X_2=\{z_1,\dots,z_{p_2}\},$$
$$X_3=\{t_1,\dots,t_{p_3}\},\quad X_4=\{u_1,\dots,u_{p_4}\},$$
where $p_0+p_1+p_2+p_3+p_4=n+3$. Then any $n$-dimensional fan
$\Sigma$ with the set of generators $\bigcup X_i$ and five
primitive collections $Y_i$ can be described up to a symmetry of
the pentagon by the following primitive relations with nonnegative
integral coefficients $c_2,\dots,c_{p_2},b_1,\dots,b_{p_3}$:
$$v_1+\dots+v_{p_0}+y_1+\dots+y_{p_1}-c_2z_2-\dots-c_{p_2}z_{p_2}-(b_1+1)t_1-\dots-(b_{p_3}+1)t_{p_3}=0,$$
$$y_1+\dots+y_{p_1}+z_1+\dots+z_{p_2}-u_1-\dots-u_{p_4}=0,$$
$$z_1+\dots+z_{p_2}+t_1+\dots+t_{p_3}=0,$$
$$t_1+\dots+t_{p_3}+u_1+\dots+u_{p_4}-y_1-\dots-y_{p_1}=0,$$
$$u_1+\dots+u_{p_4}+v_1+\dots+v_{p_0}-c_2z_2-\dots-c_{p_2}z_{p_2}-b_1t_1-\dots-b_{p_3}t_{p_3}=0.$$
\ktw
In this case we may assume that
$$v_1,\dots, v_{p_0},y_2,\dots, y_{p_1},z_2,\dots, y_{p_2},t_1,\dots, t_{p_3},u_2,\dots, u_{p_4}$$
form a basis of the lattice $N$. The other vectors are given by
\begin{align}\label{rownania}
z_1=&-z_2-\dots-z_{p_2}-t_1-\dots-t_{p_3}\notag \\
y_1=&-y_2-\dots-y_{p_1}-z_1-\dots-z_{p_2}+u_1+\dots+u_{p_4}\\
u_1=&-u_2-\dots-u_{p_4}-v_1-\dots-v_{p_0}+c_2z_2+\dots+c_{p_2}z_{p_2}\notag\\&+b_1t_1+\dots+b_{p_3}t_{p_3}\notag \\
\notag
\end{align}

\section{First results and methods}\label{sFirstResMeth}

\subsection{Bondal's construction and Thomsen's algorithm}\label{ssBondThom}

We start this section by recalling Thomsen's \cite {thom} algorithm for computing the summands of the push forward of a line bundle by a Frobenius morphism. We do this because of two reasons.

First is that Thomsen in his paper assumes finite characteristic of the ground field and uses absolute Frobenius morphism. We claim that the arguments used apply also in case of geometric Frobenius morphism and characteristic zero.

Moreover by recalling all methods we are able to show that the results of Thomsen coincide with the results stated by Bondal in \cite{ober}. Combining these both methods enables us to deduce  some interesting facts about toric varieties.

Most of the results of this section are due to Bondal and Thomsen. We use the notation from \cite{thom}. Let $\Sigma\subset N$ be a fan such that the toric variety $X=X(\Sigma)$ is smooth. Let us denote by $\sigma_i\in\Sigma$ the cones of our fan and by $T$ the torus of our variety. If we fix a basis $(e_1,\dots, e_n)$ of the lattice $N$, then of course $T=\Spec R$, where $R=k[X_{e_1^*}^{\pm 1},\dots,X_{e_n^*}^{\pm 1}]$.

In characteristic $p$ we have got two $p$-th Frobenius morphisms $F:X\rightarrow X$. One of them is the absolute Frobenius morphism given as an identity on the underlying topological space and a $p$-th power on sheaves. Notice that on the torus it is given by a map $R\rightarrow R$ that is simply a $p$-th power map, hence it is not a morphism of $k$ algebras (it is not an identity on $k$).

The other morphism is called the geometric Frobenius morphism and can be defined in any characteristic. Let us fix an integer $m$. Consider a morphism of tori $T\rightarrow T$ that associates $t^m$ to a point $t$. This is a morphism of schemes over $k$ that can be extended to the $m$-th geometric Frobenius morphism $F:X\rightarrow X$. What is important is that both of these morphisms can be considered as endomorphisms of open affine subsets associated to cones of $\Sigma$. We claim that in both cases the Thomsen's algorithm works.

We begin by recalling the algorithm from \cite{thom}.
Let $v_{i1},\dots,v_{id_i}$ be the ray generators of the $d_i$ dimensional cone $\sigma_i$. As the variety was assumed to be smooth we may extend this set to a basis of $N$. Let $A_i$ be a square matrix whose rows are vectors $v_{ij}$ in the fixed basis of $N$. Let $B_i=A_i^{-1}$ and let $w_{ij}$ be the $j$-th column of $B_i$. Of course the columns of $B_i$ are ray generators (extended to a basis) of the dual cone $\sigma_i^*\subset M=N^*$.

Let us remind that $X(\Sigma)$ is covered by affine open subsets $U_{\sigma_i}=\Spec R_i$, where $R_i=k[X^{w_{i1}},\dots,X^{w_{id_i}},X^{\pm w_{id_{i}+1}},\dots, X^{\pm w_{in}}]$. Here we use the notation $X^v=X_{e_1^*}^{v_1}\cdot\dots\cdot X_{e_n^*}^{v_n}$. Let also $X_{ij}=X^{w_{ij}}$. In this way the monomials $X_{i1},\dots,X_{in}$ should be considered as coordinates on the affine subset $U_{\sigma_i}$, so we are able to think about monomials on $U_{\sigma_i}$ as vectors: a vector $v$ corresponds to the monomial $X_i^v$. Of course all of these affine subsets contain $T$, that corresponds to the inclusions $R_i\subset R$.

Using the results of \cite{fult} we know that $U_{\sigma_i}\cap U_{\sigma_j}=U_{\sigma_i\cap\sigma_j}$ and this is a principal open subset of $U_{\sigma_i}$. This means that there is a monomial $M_{ij}$ such that $U_{\sigma_i\cap\sigma_j}=\Spec ((R_i)_{M_{ij}})$.

We are interested in Picard divisors. A $T$ invariant Picard divisor is given by a compatible collection $\{(U_{\sigma_i},X_i^{u_i})\}_{\sigma_i\in\Sigma}$. Compatible means that the quotient of any two functions in the collection is invertible on the intersection of domains. This motivates the definition:
$$I_{ij}=\{v:X_i^v\text{ is invertible in }(R_i)_{M_{ij}}\}.$$
Given a monomial $X_i^v$, if we want to know how it looks in coordinates $X_{e_1^*},\dots,X_{e_n^*}$ (obviously from the definition of $X_i$) we just have to multiply $v$ by $B_i$: $X_i^v=X^{B_iv}$. We see that $X_i^v=X_j^{B_j^{-1}B_i}$. That is why we define $C_{ij}=B_j^{-1}B_i$ and we think of $C_{ij}$ as  the matrices that translate the monomials in coordinates of one affine piece to another.

Now the compatibility in the definition of a Cartier divisor simply is equivalent to the condition $u_j-C_{ij}u_i\in I_{ji}$. We define $u_{ij}=u_j-C_{ij}u_i$ and think about them as transition maps. Of course a divisor is principal if and only if $u_{ij}=0$ for all $i,j$ (vector equal to $0$ corresponds to a constant function equal to $1$).

Let $P_m=\{v=(v_1,\dots,v_n):0\leq v_i<m\}$. Later we will see that this set has got a description in terms of characters of the kernel of the Frobenius map between tori.

Using simple algebra Thomsen proves that the following functions are well defined (the only think to prove is that the image of $h$ is in $I_{ji}$):

Let us fix $w\in I_{ji}$ and a positive integer $m$. We define the functions
$$h_{ijm}^w:P_m\rightarrow I_{ji}$$
$$r_{ijm}^w:P_m\rightarrow P_m,$$
for any $v\in P_m$ by the equation
$$C_{ij}v+w=mh_{ijm}^w(v)+r_{ijm}^w(v).$$
This is a simple division by $m$ with the rest. Moreover $r_{ijm}^w$ is bijective.

Now if we have any $v\in P_m$, a $T$-Cartier divisor $D=\{(U_{\sigma_i},X_i^{u_i})\}_{\sigma_i\in\Sigma}$ and a fixed $\sigma_l\in\Sigma$ then Thomsen defines $t_i=h_{lim}^{u_{li}}(v)$. He proves that the collection $\{(U_{\sigma_i},X_i^{t_i})\}_{\sigma_i\in\Sigma}$ is a $T$-Cartier divisor $D_v$. This is of course independent on the representation of $D$ up to linear equivalence. The choice of $l$ corresponds to "normalizing" the representation of $D$ on the affine subset $U_{\sigma_l}$. Although the definition of $D_v$ may depend on $l$, the vector bundle $\oplus_{v\in P_m}\o(D_v)$ is independent on $l$. Moreover Thomsen proves that in case of $p$-th absolute Frobenius morphism and characteristic $p>0$ this vector bundle is a push forward of the line bundle $\o(D)$. The proof uses only the fact that the Frobenius morphism can be considered as a morphism of affine pieces $U_{\sigma_i}$, so can be extended to the case of geometric Frobenius morphism and arbitrary characteristic. One only has to notice that the basis of free modules obtained by Thomsen in \cite[Section 5, Theorem 1]{thom} are exactly the same in all cases.

Now let us remind that there is an exact sequence \ref{PPic}:
$$0\rightarrow M\rightarrow D_T\rightarrow Pic\rightarrow 0,$$
where $D_T$ are $T$ invariant divisors. Let $(g_j)$ be the collection of ray generators of the fan $\Sigma$ and $D_{g_j}$ a divisor associated to the ray generator $g_j$. The morphism from $M$ to $D_T$ is given by $v\rightarrow\sum_j v(g_j)D_{g_j}$. Such a map may be extended to a map from $M_{\r}=M\otimes_{\z}\r$ by $f:v\rightarrow\sum_j [v(g_j)]D_{g_j}$. Notice that this is no longer a morphism, however if $a\in M$ and $b\in M_\r$, then $f(a+b)=f(a)+f(b)$. We obtain a map $\textbf{T}:=\frac{M_\r}{M}\rightarrow Pic$, where $\textbf{T}$ is a real torus (do not confuse with $T$). We also fix the notation for an $\r$-divisor $D=\sum_j a_jD_{g_j}$:
$$[D]:=\sum_j[a_j]D_{g_j}$$

Let $G$ be the kernel of the $m$-th geometric Frobenius morphism between the tori $T$. By acting with the functor $\Hom(\cdot,\c^*)$ we obtain an exact sequence:
$$0\rightarrow M\rightarrow M\rightarrow G^*\simeq\frac{M}{mM}\rightarrow 0.$$
We also have a morphism:
$$\frac{1}{m}:G^*\simeq\frac{M}{mM}\rightarrow\textbf T,$$
that simply divides the coordinates by $m$. By composing it with the morphism from $\textbf T\rightarrow Pic$ we get a morphism from $G^*$ to $Pic$. It can be also described as follows:

We fix $\chi\in G^*$ and arbitrarily lift it to an element $\chi_M\in M$. Now we use the morphism $M\rightarrow Div_T$ to obtain a $T$ invariant principal divisor $D_\chi$. The image of $\chi$ in $Pic$ is simply equal to $[\frac{D_{\chi}}{m}]$. Of course for different lifts of $\chi$ to $M$ we get linearly equivalent divisors. Now we prove one of the results stated by Bondal in \cite{ober}:
\ob\label{wzornapush}
Let $L=\o(D)$ by any line bundle on a smooth toric variety $X$. The push forward $F_*(\o(D))$ is equal to $\oplus_{\chi\in G^*}\o([\frac{D+D_{\chi}}{m}])$.
\kob
\uwa
The characters of $G$ play the role of $v\in P_m$ in Thomsen's algorithm.
Notice also that it is not clear that $\oplus_{\chi\in G^*}\o([\frac{D+D_{\chi}}{m}]$ is independent on the representation of $L$ by $D$. If we prove that this is equal to the push forward then this fact will follow, but in the proof we have to take any representation of $L$ and we cannot change $D$ with a linearly equivalent divisor.
\kuwa
\dow
Let $D=\{(U_{\sigma_i},X_i^{u_i})\}$ and let us fix $\chi\in G^*$. We have to prove that $\o([\frac{D+D_{\chi}}{m}])$ is one of $\o(D_v)$ for $v\in P_m$ and that this correspondence is one to one over all $\chi\in G^*$. We already know that $[\frac{D_\chi}{m}]$ is independent on the choice of the lift of $\chi$, so we may take such a lift, that $v=\chi_M+u_l$ is in the $P_m$. Here $l$ is an index of a cone, but we may assume that its ray generators form a standard basis of $N$, so $A_l=Id$. Of course such a matching between $\chi\in G^*$ and $v\in P_m$ is bijective.

Now let us compare the coefficients of $[\frac{D+D_{\chi}}{m}]$ and $D_v$. We fix a ray generator $r=(r_1,\dots,r_n)\in\sigma_j$. Let $k$ be such that this ray generator is the $k$-th row of matrix $A_j$. We compare coefficients of $D_r$. Let $\chi_M=(a_1,\dots,a_n)$. We see that:
$$\left[\frac{D+D_{\chi}}{m}\right] =\dots+\left[ \frac{(u_j)_k+\sum_{w=1}^na_wr_w}{m}\right] D_r+\dots.$$
Here of course $(u_j)_k$ is not a transition map $u_{jk}$, but the $k$-th entry of vector $u_j$ that is of course the coefficient of $D_r$ of the divisor $D$. Now from Thomsen's algorithm described above we know that
$$C_{lj}(\chi+u_l)+u_{lj}=mt_j+r,$$
where $r\in P_m$. We see that
$$t_j=\left[ \frac{C_{lj}(\chi+u_l)+u_{lj}}{m}\right] .$$
Now $A_l=Id$ and from the definition of $u_{lj}$ we have $C_{lj}u_l+u_{lj}=u_j$, so:
$$t_j=\left[\frac{A_j\chi+u_j}{m}\right] .$$
This gives us:
$$D_v=\dots+\left[ \frac{\sum_{w=1}^na_wr_w+(u_j)_k}{m}\right] D_r+\dots$$
what completes the proof.
\kdow

From \cite{ober} we know that the image $B$ of $\textbf T$ in $\Pic$ is a full collection of line bundles. Of course $B$ is a finite set (the coefficients of divisors associated to ray generators are bounded). Moreover the image of rational points of $\textbf T$ contains the whole image of $\textbf T$ (a set of equalities and inequalities with rational coefficients has got a solution in $\r$ if and only if it has got a solution in $\q$). This means that for sufficiently large $m$ the split of the push forward of the trivial bundle by the $m$-th Frobenius morphism coincides with the image of $\textbf T$ and hence is full.

Let us now consider an example of $\p^2$. Let $v_1$, $v_2$ and $v_3=-v_1-v_2$ be the ray generators of the fan. We fix a basis $(v_1,v_2)$ of $N$. The image of the torus $\textbf T$ is equal to the set of all divisors of the form $[a]D_{v_1}+[b]D_{v_2}+[-a-b]D_{v_3}$ for $0\leq a,b<1$. We see that the image of the torus $\textbf T$ is $\o,\o(-1),\o(-2)$. This is a full collection. Notice however that it is not true that if we have a line bundle $L$ then there exists an integer $m_0$ such that the push forward of $L$ by the $m$-th Frobenius morphism for $m>m_0$ is a direct sum of line bundles from $B$. For example the push forward of $\o(-3)$ always contains in the split $\o(-3)$ that is not an element of $B$. However, as we will see only minor differences from the set $B$ are possible.
\dfi
Let us fix a natural bijection between points of $\textbf T$ and elements of $M_\r$ with entries from $[0,1)$ in some fixed basis. Now each element of $B$ has got a natural representant in $Div_T$ as sum of $D_{g_j}$ with integer coefficients. Let $B_0\subset Div_T$ be the set of these representatives. We define the set $B'$ as the set of all divisors $D$ in $\Pic$ for which there exists an element in $b\in B_0$, such that there exists a representation of $D$ whose coefficients differ by at most one from the coefficients of $b$.
\kdfi

In other words we take (some fixed) representations of all elements of $B$, we take all other representations whose coefficients differ by at most one and we take the image in $\Pic$ to obtain $B'$.

Let us look once more at the example of $\p^2$. With previous notation $B$ is equal to $0$, $-D_{v_3}$, $-2D_{v_3}$. The set $B'$ would be equal to $\pm D_{v_1}\pm D_{v_2}\pm D_{v_3}$, $\pm D_{v_1}\pm D_{v_2}\pm D_{v_3}-D_{v_3}$, $\pm D_{v_1}\pm D_{v_2}\pm D_{v_3}-2D_{v_3}$. This gives us $\o(3)$,$\o(2)$,$\o(1)$, $\o$, $\o(-1)$, $\o(-2)$, $\o(-3)$, $\o(-4)$, $\o(-5)$,.
\ob
For any smooth toric variety and any line bundle there exists an integer $m_0$ such that the push forward by the $m$-th Frobenius morphism for any $m>m_0$ splits into the line bundles form $B'$.
\kob
\dow
From \ref{wzornapush} we know that the line bundles from the split are of the form $[\frac{D}{m}+\frac{D_\chi}{m}]$, where $L=\o(D)$ is a fixed representation of $L$. Of course for sufficiently large $m$ all coefficients of $\frac{D}{m}$ belong to the interval $(-1,1)$, so the coefficients of $[\frac{D}{m}+\frac{D_\chi}{m}]$ differ by at most one from the coefficients of $[\frac{D_\chi}{m}]$ that is in $B$, so in fact $[\frac{D}{m}+\frac{D_\chi}{m}]\in B'$.
\kdow
This combined with the result of Thomsen \cite{thom} that the push forward and the line bundle are isomorphic as sheaves or abelian groups gives us the following result:
\wn\label{dz}
There exists a finite set, namely $B'$, such that each line bundle is isomorphic as a sheaf of abelian groups to a direct sum of line bundles from $B'$. In particular their cohomologies agree.
\kwn
\subsection{Techniques of counting homology}
Our aim will be to describe line bundles on toric varieties with vanishing
higher cohomologies, that we call acyclic. Later, we will use this
characterization to check if $Ext^i(L,M)=H^i(L^\vee\otimes M)$ is equal to zero
for $i>0$. We start with general remarks on cohomology of line bundles on smooth, complete toric varieties.

Let $\Sigma$ be a fan in $N=\mathbb Z^n$ with rays $x_1,...,x_m$
and let $\mathbb P_{\Sigma}$ denote the variety constructed from
the fan $\Sigma$. For $I\subset\{1,\dots,m\}$ let $C_I$ be a
simplicial complex generated by sets $J\subset I$ such that
$\{x_i:i\in J\}$ generate a cone in $\Sigma$. For
$r=(r_i:i=1,\dots,m)$ let us define $Supp(r):=C_{\{i:\;r_i\geq
0\}}$.

The proof of the following well known fact can be found in the paper \cite{bohu}:
\ob\label{cohogen} The cohomology $H^j(\mathbb P_{\Sigma},L)$ is
isomorphic to the direct sum over all $r=(r_i:i=1,\dots,m)$ such
that $\o(\sum_{i=1}^m r_iD_{x_i})\cong L$ of the $(n-j)$-th
reduced homology of the simplicial complex $Supp(r)$.\kob

\dfi We call a line bundle $L$ on $\mathbb P_{\Sigma}$ acyclic if
$H^i(\mathbb P_{\Sigma},L)=0$ for all $i\geq 1$.\kdfi

\dfi For a fixed fan $\Sigma$ we call a proper subset $I$ of $\{1,\dots,m\}$ a forbidden
set if the simplicial complex $C_I$ has nontrivial reduced
homology.\kdfi

From Proposition \ref{cohogen} we have the following
characterization of acyclic line bundles:

\ob\label{acycliclb} A line bundle $L$ on $\mathbb P_{\Sigma}$ is acyclic if it is
not isomorphic to any of the following line bundles
$$\o(\sum_{i\in I}r_iD_{x_i}-\sum_{i\not\in I}(1+r_i)D_{x_i})$$
where $r_i\geq 0$ and $I$ is a proper forbidden subset of
$\{1,\dots,m\}$. \kob

Hence to determine which bundles on $\mathbb P_{\Sigma}$ are
acyclic it is enough to know which sets $I$ are forbidden.

In our case $C_I=\{J\subset I: \widehat{Y_i}:=\{j:x_j\in Y_i\}\nsubseteq J$ for
$i=1,\dots,5\}$, since $Y_i$ are primitive collections. We call
sets $\widehat{Y_i}$ also primitive collections. The only difference
between sets $\widehat{Y_i}$ and $Y_i$ is that the first one is the set of
indices of rays in the second one, so in fact they could be even
identified.

In case of simplicial complex $S$ on the set of vertices $V$ we also define a primitive collection as a minimal subset of vertices that do not form a simplex. Complex $S$ is determined by its primitive collections, namely it contains simplexes (subsets of $V$) that contain none of primitive collections.

We describe very powerful method of counting homologies of
simplicial complexes which are given by their primitive
collections (as in our case). We use the result of
Mrozek and Batko \cite{mroz}:

\lem\label{wyrwac} Let $X$ be a simplicial complex and let $Z$ be
a cycle in the chain complex whose boundary $B$ is exactly one simplex.
Then we can remove the pair $(Z,B)$ from the chain complex without
changing the homology.\klem

\dfi Let $X$ be a simplicial complex defined by its set of
primitive collections $\textsc{P}$ on the set of vertices $V$. We
say that simplicial complex $X'$ on the set of vertices
$V\setminus P$ is obtained from $X$ by delating a primitive
collection $P$ if the set of primitive collections of $X'$ is
equal to the set of minimal sets in $\{Q\cap(X\setminus P):Q\in
\textsc{P}\}$.\kdfi

\lem\label{usunac} Let $X$ be a simplicial complex and suppose
that there exists an element $x$ which belongs to exactly one
primitive collection $P$. Let $m=|P|$ and let $X'$ be a simplicial
complex obtained from $X$ by delating $P$, then
$$h^i(X)=h^{i-m+1}(X').$$
\klem \dow Using Lemma \ref{wyrwac} we will be removing subsequently on
dimension reductive pairs $(Z,B)$ such that $x\in Z$. We start
from $(\{x\},\emptyset)$. One can see that in each dimension we
can take all $(Z,Z\setminus\{x\})$ for $Z$ containing $x$ as
reductive pairs. Let us consider all simplexes of $X$ that do not
contain $P\setminus\{x\}$. One can prove by induction on dimension
that we will remove all of them:

Let $D$ be a simplex. If it contains $x$, than it will be removed
as a first element of a reductive pair. If it does not, then $D\cup
\{x\}$ is also a simplex of $X$ and we will remove
$(D\cup\{x\},D)$.

We see that our simplicial complex can be reduced to a complex
with simplexes containing $P\setminus\{x\}$. Now one immediately
sees that such a complex is isomorphic to a complex $X'$ (with a
degree shifted by $|P\setminus\{x\}|=m-1$).\kdow

The same method allows us to easily compute homologies when there
are few primitive collections and many points. The idea is that we
can glue together points that are in exactly the same primitive
collections.

\dfi Let $X$ be a simplicial complex defined by its set of
primitive collections $\textsc{P}$ on the set of vertices $V$.
Suppose that there exist two points $x,y\in X$ such that they
belong to the same primitive collections. We say that a simplicial
complex $X'$ on the set of vertices $V\setminus\{y\}$ is obtained
from $X$ by gluing points $x$ and $y$ if the set of primitive
collections of $X'$ is equal $\{Q\setminus\{y\}:Q\in\textsc{P}\}$.
We can think of it like $x$ was in fact two points $x,y$.\kdfi

\ob\label{red} Let $X$ be a simplicial complex and suppose that
there exist two points $x,y\in X$ such that they belong to the
same primitive collections. Let $X'$ be a simplicial complex
obtained from $X$ by gluing points $x$ and $y$, then
$$h^i(X)=h^{i-1}(X').$$\kob \dow In both complexes we will be
removing reductive pairs of the form $(Z,B)$ with $x\in Z$ just as
in Lemma \ref{usunac}. In both situations all that is left are simplexes
that contain a set of a form $P\setminus\{x\}$, where $P$ is a
primitive collection containing $x$. In this situation all of
simplexes of $X$ that are left contain $y$ and they can be
identified with simplexes of $X'$ that are left, the maps are
exactly the same what finishes the proof. \kdow

\wn\label{skleic} Let $X$ be a simplicial complex on the set of
vertices $V$. Let $X'$ be a simplicial complex obtained from $X$
by gluing equivalence classes of the relation $\sim $ that identifies elements that are in exactly
the same primitive collections. Suppose $|V|-|V/\sim|=m$, then
$$h^i(X)=h^{i-m}(X').$$\kwn \dow We use \ref{red}
for pairs of points in the equivalence classes.\kdow

\wn In the situation of Lemma \ref{usunac} and Corollary
\ref{skleic} $X$ is acyclic if and only if $X'$ is acyclic.\kwn

With these tools we are ready to determine forbidden subsets. In general we have got two following Lemmas:

\lem If a nonempty subset $I$ is not a sum of
primitive collections, then it is not forbidden. \klem \dow There
exists $a\in I$ such that $a$ does not belong to any primitive
collection which is contained in $I$. Using Lemma \ref{wyrwac} we can remove subsequently on
dimension reductive pairs $(Z,B)$ such that $a\in Z$. We start
from $(\{a\},\emptyset)$. One can see that in this way we remove all of simplexes and as a consequence the chain complex is exact.\kdow

\lem\label{jedna} A primitive collection is a forbidden subset. \klem \dow Using Lemma \ref{usunac} we can remove this primitive collection and get a complex consisting of the empty set only that has nontrivial reduced homologies.

This can be also seen from the fact that the considered complex topologically is a sphere.
\kdow

The following Lemmas apply to the case when the Picard number is three and we have five primitive collections as in Batyrev's classification. Let us remind that primitive collections of simplicial compex in this case are $\widehat{Y_i}:=\{j:x_j\in Y_i\}$, for our convenience we define also $\widehat{X_i}:=\{j:x_j\in X_i\}$.

\lem A sum of two consecutive primitive collections is a forbidden
subset. \klem \dow
Using Lemma \ref{usunac} we remove one primitive collection and get a situation of Lemma \ref{jedna}.
\kdow

\lem A sum of three consecutive primitive collections $\widehat{Y_{i}}$, $\widehat{Y_{i+1}}$, $\widehat{Y_{i+2}}$ is not a
forbidden subset. \klem \dow First we can remove primitive collection
$\widehat{Y_{i}}$. The image of $\widehat{Y_{i+2}}$ contains the image of $\widehat{Y_{i+1}}$, so in fact
we are left with just one primitive collection $P$ which is an image of $\widehat{Y_{i+1}}$. We can remove $P$ and obtain a nonempty full simplicial complex which is known to have trivial homologies.\kdow

Above Lemmas match together to the following

\tw\label{forbiddensets} The only forbidden subsets are primitive collections, their
complements and the empty set.\ktw

This gives us that in our situation

\wn\label{acyclic} A line bundle $L$ is acyclic if and only if it is
not isomorphic to any of the following line bundles
$$\o(\alpha_1^1D_{v_1}+\cdots+\alpha_2^1D_{y_1}+\cdots+\alpha_3^1D_{z_1}+\cdots+\alpha_4^1D_{t_1}+\cdots+\alpha_5^1D_{u_1}+\cdots)$$
where exactly $2,3$ or $5$ consecutive $\alpha_i:=(\alpha_i^1,\cdots,\alpha_i^{p_i})$ are all less or equal to $-1$ and the rest is nonnegative.\kwn \dow  It is an immediate consequence of Proposition \ref{acycliclb} and Theorem \ref{forbiddensets}\kdow

\wn\label{acyclic2} If all of the coefficients $b$ and $c$ are zero in the primitive relations from Theorem \ref{charBatyrev} then a line bundle $L$ is acyclic if and only if it is
not isomorphic to any of the following line bundles
$$\o(\alpha_1D_{v}+\alpha_2D_{y}+\alpha_3D_{z}+\alpha_4D_{t}+\alpha_5D_{u})$$
where exactly $2,3$ or $5$ consecutive $\alpha_i$ are negative and
if $\alpha_i<0$ then $\alpha_i\leq-|X_i|$.\kwn \dow Since all
divisors corresponding to elements of the set $X_i$ are linearly
equivalent we match them together and as a consequence $\alpha_i$
is the sum of all of their coefficients. \kdow

\section{Main theorem}\label{yes}

This section contains the main, new result of this work. We give an
explicit construction of a full, strongly exceptional collection of
line bundles in the derived category $D^b(X)$ for a large family of
smooth, complete toric varieties $X$ with Picard number three.
Namely for varieties $X$ whose sets $X_1,X_3$ and $X_4$ from
Batyrev's classification \ref{charBatyrev} have only one element. We will use
results from Section \ref{sFirstResMeth}.

\subsection{Our setting}\label{ssour}

In this subsection we establish a family of varieties which we
consider in this section, we also fix notation.

From now on for the whole Section let $X$ be smooth, complete
toric variety with Picard number three, which using the notation
from Theorem \ref{charBatyrev} has $|X_1|=|X_3|=|X_4|=1$.

Let $r=|X_2|$. Then of course $|X_0|=n-r$. We allow arbitrary
nonnegative integer parameters $b:=b_1$, $c_2,\dots,c_r$. This
family generalizes one considered in \cite{dlm} (there, the case
$r=1$ was considered) and \cite{lcmr} (there the case
$b=c_1=\dots=c_r=0$ was considered).

\uwa
A variety of this type is Fano iff
$$n-r> \sum_{i=2}^r c_r+b.$$
In what follows we do not restrict to the Fano case.
\kuwa

First let us write what are the coordinates of the ray generators
in the considered situation:

\begin{align}\label{raywbazie}
&v_1=e_1, v_2=e_2,\dots, v_{n-r}=e_{n-r}\notag\\
&y=-e_1-\dots-e_{n-r}+c_2e_{n-r+2}+\dots+c_re_n-(b+1)(e_{n-r+1}+\dots+e_n)\notag\\
&z_1=e_{n-r+1},\dots,z_r=e_n\\
&t=-e_{n-r+1}-\dots-e_n\notag\\
&u=-e_1-\dots-e_{n-r}+c_2e_{n-r+2}+\dots+c_re_n-b(e_{n-r+1}+\dots+e_n)\notag\\
\notag\end{align}

Let $D_w$ be the divisor associated to the ray generator $w$. One
can easily see that the divisors $D_{v_1},\dots,D_{v_{n-r}}$ are
all linearly equivalent. Let $D_v$ be any their representant in
the Picard group. The other equivalence relations that generate
all the relations in the Picard group are:

\begin{align}\label{relacje}
D_v&\simeq D_u+D_y\notag\\
D_{z_1}&\simeq D_t+bD_u+(b+1)D_y\\
D_{z_i}&\simeq D_t+(b-c_i)D_u+(b-c_i+1)D_y\quad 2\leq i\leq
r\notag\\\notag
\end{align}

From these relations we can easily deduce:

\ob\label{Picard} The Picard group of the variety $X$ is
isomorphic to $\z^3$ and is generated by $D_t,D_y,D_v$.\kob

We introduce two sets of
divisors. We claim that these sets can be ordered in
such a way that line bundles corresponding to divisors from these
sets form a strongly exceptional collection.

\begin{align}\label{kolekcja}
Col_1=\{&-sD_t-sD_y+(-(n-r)-bs+q)D_v:\notag\\&0\leq s\leq r,0\leq q\leq n-r\}\\
Col_2=\{&-sD_t-(s-1)D_y+(-(n-r)-bs+q)D_v:\notag\\&1\leq s\leq r,0\leq q\leq n-r-1\}\notag\\\notag
\end{align}

\dfi
Let $Col=Col_1\cup Col_2$.
\kdfi

\uwa\label{ilewybranych} Let us notice that $|Col_1|=(r+1)(n-r+1)$
and $|Col_2|=r(n-r)$, so $|Col|=2rn-2r^2+n+1$.\kuwa

We calculate the number of maximal cones in the fan defining the
variety $X$. In order to obtain a maximal cone we have to choose
$n$ ray generators that do not contain a primitive collection.
This is equivalent to removing three ray generators in such a way
that the rest do not contain a primitive collection. First let us
notice that we can remove at most one element from each group
$X_i$ because otherwise the rest would contain a primitive
collection. We have the following possibilities:

1) We remove one element from $X_0$ and $X_2$. Then we have to
remove one element from $X_3$ or $X_4$. We have got $2(n-r)r$ such
possibilities.

2) We remove one element from $X_0$ and none from $X_2$. We have
got $n-r$ such possibilities.

3) We remove one element from $X_2$ and none from $X_0$. We have
got $r$ such possibilities.

4) We do not remove any elements from $X_0$ and from $X_2$. We
have got $1$ such possibility.

All together we see that we have $2rn-2r^2+n+1$ maximal cones.
From the general theory we know that the rank of the Grothendieck
group is the same. Let us notice that from Remark
\ref{ilewybranych} our set $Col$ is of the same number of
elements.

\subsection{Acyclicity of differences of line bundles from
$Col$}\label{ssAcyc}

In this Subsection we order the set $Col$ and prove that line
bundles corresponding to divisors from $Col$ form a strongly
exceptional collection.

Let us first check that $\Ext^i_{\o_X}(\o(D_1),\o(D_2))=0$ for any
divisors $D_1,D_2$ from the set $Col$ and for any $i>0$. We know
that
$$\Ext^i_{\o_X}(\o(D_1),\o(D_2))=H^i(\o(D_1)^\vee\otimes\o(D_2))=H^i(\o(D_2-D_1)).$$
This means that we have to show  that all line bundles associated
to differences of divisors from $Col$ are acyclic.

\dfi Let $Diff$ be the set of all divisors of the form $D_1-D_2$,
where $D_1,D_2\in Col$. \kdfi

\ob\label{roznice} The set $Diff$ is the sum of sets
$Diff_1,Diff_2,Diff_3$, where:
$$Diff_1=\{sD_t+sD_y+(bs+q)D_v:$$
$$-r\leq s\leq r,r-n\leq q\leq n-r\}$$
$$Diff_2=\{sD_t+(s-1)D_y+(bs+q)D_v:$$
$$-r+1\leq s\leq r,r-n+1\leq q\leq n-r\}$$
$$Diff_3=\{sD_t+(s+1)D_y+(bs+q)D_v:$$
$$-r\leq s\leq r-1,r-n\leq q\leq n-r-1\}$$
\kob \dow The set $Diff_1$ is equal to the set of all possible
differences of two divisors from $Col_1$ and this set contains all
possible differences of two divisors from $Col_2$. The set
$Diff_2$ is the set of all possible differences of the form
$D_1-D_2$, where $D_1\in Col_1,D_2\in Col_2$. The set $Diff_3$ is
equal to $-Diff_2$ and so it is equal to the set of all differences
of the form $D_2-D_1$, where $D_1\in Col_1,D_2\in Col_2$. This are
of course all possible differences of two elements form $Col$.
\kdow

From the Corollary \ref{acyclic} we know that it is enough to
prove that elements of $Diff$ are not of the form
\begin{align*}
\alpha_1D_v+\alpha_2D_y+\alpha_3^1D_{z_1}+\alpha_3^2D_{z_2}+\dots+\alpha_3^r D_{z_r}+\alpha_4 D_t+\alpha_5 D_u,
\end{align*}
where exactly two, three or five consecutive $\alpha_i$'s are negative (we
call a number positive when it is nonnegative and consider only
two signs positive and negative) and:

1) if $\alpha_1<0$, then $\alpha_1\leq -(n-r)$ ($\alpha_1$ is in fact
sum of all the coefficients of $D_{v_i}$, which have to be of the
same sign),

2) if any $\alpha_3^i<0$ then $\alpha_3^j<0$ (all parameters
$\alpha_3^j$ are treated as one group and have the same sign).

From now on we assume that these conditions on $\alpha_i$'s are
satisfied.

Using the relations \ref{relacje} we obtain:
\begin{align}\label{forma}
&\alpha_1D_v+\alpha_2D_y+\alpha_3^1D_{z_1}+\alpha_3^2D_{z_2}+\dots+\alpha_3^r D_{z_r}+\alpha_4 D_t+\alpha_5 D_u=\notag\\
&(\alpha_4+\sum_{j=1}^r\alpha_3^j)D_t+(\alpha_2-\alpha_5+\sum_{j=1}^r\alpha_3^j)D_y+\notag\\
&(\alpha_1+b\alpha_3^1+\sum_{j=2}^r(b-c_j)\alpha_3^j+\alpha_5)D_v\\\notag
\end{align}

\lem If the elements $\alpha_3^j$ are negative then the divisors
form $Diff$ are not of the form \ref{forma}. \klem \dow If
$\alpha_4$ was negative, then the coefficient of $D_t$ would be
less then or equal to $-r-1$ and none of the divisors from $Diff$
has got such a coefficient, so $\alpha_4$ has to be positive.
Since $\alpha_3$ is negative and $\alpha_4$ is positive, then
$\alpha_2$ has to be negative and $\alpha_5$ has to be positive.
This means that the coefficient of $D_y$ is less then or equal to
$-r-1$. The divisors from $Diff$ are not of this form. \kdow

From now on we may assume that $\alpha_3$ is positive.

\lem\label{jeden} The divisors from $Diff_1$ are not of the form
\ref{forma}. \klem \dow Suppose that a divisor from $Diff_1$ can
be written in a form \ref{forma}. We have:
$$\alpha_4+\sum_{j=1}^r\alpha_3^j=\alpha_2-\alpha_5+\sum_{j=1}^r\alpha_3^j,$$
so $\alpha_4+\alpha_5=\alpha_2$. But $\alpha_2,\alpha_4$ and
$\alpha_5$ cannot be of the same sign, so $\alpha_4$ and
$\alpha_5$ have to have different signs. As $\alpha_3$ was
positive we see that $\alpha_4$ is positive, so $\alpha_5$ and
$\alpha_1$ are negative. Let us notice that:
$$\alpha_1+b\alpha_3^1+(\sum_{j=2}^r(b-c_j)\alpha_3^j)+\alpha_5\leq$$
$$-n+r+b(\sum_{j=1}^r\alpha_3^j)-1\leq$$
$$-n+r-1+b(\alpha_4+\sum_{j=1}^r\alpha_3^j)$$
This shows precisely that the coefficient of $D_v$ is less then or
equal to $-n+r-1$ plus $b$ times the coefficient of $D_t$. Let $s$
be the coefficient of $D_t$. From the definition of $Diff_1$ the
coefficient of $D_v$ is at least $-n+r+bs$. This
gives us a contradiction. \kdow

\lem\label{dwa} The divisors from $Diff_3$ are not of the form
\ref{forma}. \klem \dow Suppose that a divisor from $Diff_3$ can
be written in a form \ref{forma}. We have:
$$\alpha_4+\sum_{j=1}^r\alpha_3^j=\alpha_2-\alpha_5-1+\sum_{j=1}^r\alpha_3^j,$$
so $\alpha_4+\alpha_5=\alpha_2-1$. The rest of the proof is
identical to the proof of Lemma \ref{jeden}. \kdow \lem The
divisors from $Diff_2$ are not of the form \ref{forma}. \klem \dow
Suppose that a divisor from $Diff_2$ can be written in a form
\ref{forma}. We have:
$$\alpha_4+\sum_{j=1}^r\alpha_3^j=\alpha_2-\alpha_5+1+\sum_{j=1}^r\alpha_3^j,$$
so $\alpha_4+\alpha_5=\alpha_2+1$. But $\alpha_2,\alpha_4$ and
$\alpha_5$ cannot be of the same sign, so we have to possible
cases:

1) The coefficients $\alpha_4$ and $\alpha_5$ have different
signs. In this case the proof is the same as in Lemmas \ref{jeden}
and \ref{dwa}.

2) We have $\alpha_4=\alpha_5=0$ and $\alpha_2=-1$. In this case
$\alpha_1$ has to be negative, because $\alpha_3$ was positive.
Let $s=\alpha_4+\sum_{j=1}^r\alpha_3^j$ be the coefficient of
$D_t$. We have:
$$\alpha_1+b\alpha_3^1+\sum_{j=2}^r(b-c_j)\alpha_3^j+\alpha_5\leq-n+r+bs,$$
so the coefficient of $D_v$ is less then or equal to $-n+r+bs$.
But from the definition of $Diff_2$ we know that the coefficient
of $D_v$ is at least $bs+r-n+1$ what gives us a contradiction.
\kdow

Now we only have to order
the line bundles corresponding to divisors from $Col$ in such a way that
$$0=\Ext^0_{\o_X}(\o(D_1),\o(D_2))=H^0(\o(D_1)^\vee\otimes\o(D_2))=H^0(\o(D_2-D_1)).$$
for any divisors $D_1>D_2$.

Let us define the order by:
$L_{s,q}<L'_{s,q}<L_{s,q+1},\;L_{s+1,q_1}<L_{s,q_2}$ where
$$L_{s,q}=\o(-sD_t-sD_y+(q-bs-(n-r))D_v)$$
for $s=0,\dots,r$ and $q=0,\dots,n-r$ and
$$L'_{s,q}=\o(-sD_t-(s-1)D_y+(q-bs-(n-r))D_v)$$
for $s=1,\dots,r-1$ and $q=0,\dots,n-r-1$. It is easy to
see that zero cohomology of appropriate difference vanish.

\subsection{Generating the derived category}\label{ssGenDC}

We prove that the strongly exceptional collection from Subsection
\ref{ssour} is also full. We show that it generates all line bundles and due to the result of \cite{boho} it is enough. In order to show that we
need several lemmas:

\lem\label{1} Let $s$ and $k$ be any integers. Line bundles
$L_{q}=\o(-sD_t-sD_y+(k+q)D_v$ for $q=0,\dots,n-r$ and
$L'_{q}=\o(-sD_t-(s-1)D_y+(k+q)D_v)$ for $q=0,\dots,n-r-1$
generate $\o(-sD_t-(s-1)D_y+(n-r+k)D_v)$ in the derived
category.\klem \dow We consider the Koszul complex for
$\o(D_{y}),\o(D_{v_1}),\dots,\o(D_{v_{n-r}})$:
$$0\rightarrow \o(-D_y-(n-r)D_v)\rightarrow\dots\rightarrow\o(-D_v)^{n-r}\oplus\o(-D_y)\rightarrow\o\rightarrow 0.$$
By tensoring it with $\o(-sD_t-(s-1)D_y+(k+n-r)D_v)$ we obtain:
$$0\rightarrow \o(-sD_t-sD_y+kD_v)\rightarrow\dots\rightarrow\o(-sD_t-(s-1)D_y+(k+n-r-1)D_v)^{n-1}$$
$$\oplus\o(-sD_t-sD_y+(k+n-r)D_v)\rightarrow\o(-sD_t-(s-1)D_y)+(k+n-r)D_v)\rightarrow 0.$$
All sheaves that appear in this exact sequence, apart from the
last one, are exactly
$\o(-sD_t-sD_y+kD_v),\dots,\o(-sD_t-sD_y+(k+n-r)D_v),\o(-sD_t-(s-1)D_y+kD_v),\dots,\o(-sD_t-(s-1)D_y+(k+n-r-1)D_v)$,
so indeed we can generate $\o(-sD_t-(s-1)D_y+(k+n-r)D_v)$. \kdow

\lem\label{2} Let $s$ and $k$ be any integers. Line bundles
$L_{q}=\o(-sD_t-sD_y+(k+q)D_v)$ for $q=0,\dots,n-r$ and
$L'_{q}=\o(-sD_t-(s-1)D_y+(k+q)D_v)$ for $q=1,\dots,n-r$ generate
$\o(-sD_t-(s-1)D_y+kD_v)$ in the derived category.\klem \dow The
proof is similar to the last one. We deduce assertion from the
same exact sequence of sheaves.\kdow

\lem\label{3} Let $s$ and $k$ be any integers. Line bundles
$L_{q}=\o(-sD_t-sD_y+(k+q)D_v)$ for $q=1,\dots,n-r$ and
$L'_{q}=\o(-sD_t-(s-1)D_y+(k+q)D_v)$ for $q=0,\dots,n-r$ generate
$\o(-sD_t-sD_y+(n-r+k+1)D_v)$ in the derived category.\klem \dow
The proof is similar to the first one. We have to consider the
Koszul complex for line bundles
$\o(D_{u}),\o(D_{v_1}),\dots,\o(D_{v_{n-r}})$:
$$0\rightarrow \o(-D_u-(n-r)D_v)\rightarrow\dots\rightarrow\o(-D_v)^{n-r}\oplus\o(-D_u)\rightarrow\o\rightarrow 0$$
we dualize it and we tensor it with $\o(-sD_t-(s-1)D_y+kD_v)$.\kdow

\lem\label{4} Let $s$ and $k$ be any integers. Line bundles
$L_{q}=\o(-sD_t-sD_y+(k+q)D_v)$ for $q=1,\dots,n-r+1$ and
$L'_{q}=\o(-sD_t-(s-1)D_y+(k+q)D_v)$ for $q=1,\dots,n-r$ generate
$\o(-sD_t-sD_y+kD_v)$ in the derived category.\klem \dow The proof
is similar to the last one. We deduce assertion from the same
exact sequence of sheaves.\kdow

\lem\label{5} Let $s$ and $k$ be any integers. Line bundles
$L_{q}=\o(-sD_t-sD_y+(k+q)D_v)$ for $q=0,\dots,n-r$ and
$L'_{q}=\o(-sD_t-(s-1)D_y+(k+q)D_v)$ for $q=0,\dots,n-r-1$
generate in the derived category line bundles
$\o(-sD_t-sD_y+q'D_v)$ and $\o(-sD_t-(s-1)D_y+q'D_v)$ for
an arbitrary integer $q'$.\klem \dow We prove it by induction on
$|q'|$. For $q'\geq k+n-r$ we use Lemmas \ref{1} and \ref{3}, for
$q'<k$ we use Lemmas \ref{2} and \ref{4}.\kdow

\lem\label{6} Let $k$ be any integer. Line bundles
$L_{s,q}=\o(-sD_t-sD_y+qD_v)$ for $s=k,\dots,k+r$ and arbitrary
$q$ and $L'_{s,q}=\o(-sD_t-(s-1)D_y+qD_v)$ for $s=k,\dots,k+r-1$
and arbitrary $q$ generate in the derived category line bundles
$L'(k+r,q)=\o(-(k+r)D_t-(k+r-1)D_y+qD_v)$ with arbitrary $q$.\klem
\dow Consider the Koszul complex for
$\o(D_{y}),\o(D_{z_1}),\dots,\o(D_{z_{r}})$:
$$0\rightarrow \o(-D_{z_1}-(r-1)D_{z_2}-D_y)\rightarrow\dots$$
$$\dots\rightarrow\o(-D_{z_1})\oplus\o(-D_{z_2})^{r-1}\oplus\o(-D_y)\rightarrow\o\rightarrow 0.$$
After tensoring it with $\o(-(k-1)D_y+q'D_v)$ for appropriate $q'$
we get the assertion.\kdow

\lem\label{7} Let $k$ be any integer. Line bundles
$L_{s,q}=\o(-sD_t-sD_y+qD_v)$ for $s=k,\dots,k+r$ and arbitrary
$q$ and $L'_{s,q}=\o(-sD_t-(s-1)D_y+qD_v)$ for $s=k+1,\dots,k+r$
and arbitrary $q$ generate in the derived category line bundles
$L'(k,q)=\o(-kD_t-(k-1)D_y+qD_v)$ for arbitrary $q$.\klem \dow The
proof is similar to the last one. We deduce assertion from the
same exact sequence of sheaves.\kdow

\lem\label{8} Let $k$ be any integer. Line bundles
$L_{s,q}=\o(-sD_t-sD_y+qD_v)$ for $s=k+1,\dots,k+r$ and arbitrary
$q$ and $L'_{s,q}=\o(-sD_t-(s-1)D_y+qD_v)$ for $s=k+1,\dots,k+r+1$
and arbitrary $q$ generate in the derived category line bundles
$L(k,q)=\o(-kD_t-kD_y+qD_v)$ for arbitrary $q$.\klem \dow Consider
the Koszul complex for
$\o(D_{z_1}),\dots,\o(D_{z_{r}}),\o(D_{t})$:
$$0\rightarrow \o(-D_{z_1}-(r-1)D_{z_2}-D_t)\rightarrow\dots$$
$$\dots\rightarrow\o(-D_{z_1})\oplus\o(-D_{z_2})^{r-1}\oplus\o(-D_t)\rightarrow\o\rightarrow 0.$$
After tensoring it with $\o(-kD_y+q'D_v)$ for appropriate $q'$ we
get the assertion.\kdow

\lem\label{9} Let $k$ be any integer. Line bundles
$L_{s,q}=\o(-sD_t-sD_y+qD_v)$ for $s=k,\dots,k+r$ and arbitrary
$q$ and $L'_{s,q}=\o(-sD_t-(s-1)D_y+qD_v)$ for $s=k+1,\dots,k+r$
and arbitrary $q$ generate in the derived category line bundles
$L'(k+r+1,q)=\o(-(k+r+1)D_t-(k+r)D_y+qD_v)$ for arbitrary
$q$.\klem \dow The proof is similar to the last one. We deduce
assertion from the same exact sequence of sheaves.\kdow

\lem\label{10} Let $k$ be any integer. Line bundles
$L_{s,q}=\o(-sD_t-sD_y+qD_v)$ for $s=k,\dots,k+r$ and arbitrary
$q$ and $L'_{s,q}=\o(-sD_t-(s-1)D_y+qD_v)$ for $s=k,\dots,k+r-1$
and arbitrary $q$ generate in the derived category line bundles
$L(s,q)=\o(-sD_t-sD_y+qD_v)$ and $L'(s,q)=\o(-sD_t-(s-1)D_y+qD_v)$
for arbitrary $s$ and $q$.\klem \dow We prove it by induction on
$|s|$. For $s\geq k+n-r$ we use Lemmas \ref{6} and \ref{9}, for
$r<k$ we use Lemmas \ref{7} and \ref{8}.\kdow

\lem\label{11} Let $k$ be any integer. Line bundles
$\o(-sD_t-(s+k)D_y+qD_v)$ and $\o(-sD_t-(s+k+1)D_y+qD_v)$ for
arbitrary $s$ and $q$ generate in the derived category line
bundles $\o(-sD_t-(s+k+2)D_y+qD_v)$ for arbitrary $s$ and
$q$.\klem \dow Consider the Koszul complex for
$\o(D_{t}),\o(D_{u})$:
$$0\rightarrow \o(-D_{t}-D_{u})\rightarrow\o(-D_{t})\oplus\o(-D_{u})\rightarrow\o\rightarrow 0.$$
After tensoring it with $\o(-k'D_y+q')$ for appropriate $k'$ and
$q'$ we get the assertion.\kdow

\lem\label{12} Let $k$ be any integer. Line bundles
$\o(-sD_t-(s+k)D_y+qD_v)$ and $\o(-sD_t-(s+k+1)D_y+qD_v)$ for
arbitrary $s$ and $q$ generate in the derived category line
bundles $\o(-sD_t-(s+k-1)D_y+qD_v)$ for arbitrary $s$ and
$q$.\klem \dow Consider the Koszul complex for
$\o(D_{t}),\o(D_{u})$:
$$0\rightarrow \o(-D_{t}-D_{u})\rightarrow\o(-D_{t})\oplus\o(-D_{u})\rightarrow\o\rightarrow 0.$$
After tensoring it with $\o(-k'D_y+q')$ for appropriate $k'$ and
$q'$ we get the assertion.\kdow

\ob\label{13} Line bundles
$$L_{s,q}=\o(-sD_t-sD_y+(q-bs-(n-r))D_v)$$
for $s=0,\dots,r$ and $q=0,\dots,n-r$ and
$$L'_{s,q}=\o(-sD_t-(s-1)D_y+(q-bs-(n-r))D_v)$$
for $s=0,\dots,r-1$ and $q=0,\dots,n-r-1$ generate in the
derived category all line bundles.\kob \dow We use Lemmas \ref{5},
\ref{10}, \ref{11} and \ref{12}.\kdow

Summarizing, we have proved:

\tw\label{glowne} Let $X$ be a smooth, complete, $n$ dimensional
toric variety with Picard number three and the set of ray
generators $X_0\cup\dots\cup X_4$, where
$$X_0=\{v_1,\dots,v_{n-r}\},\hskip 3pt X_1=\{y\},\hskip 3pt X_2=\{z_1,\dots,z_{r}\},\hskip 3pt X_3=\{t\},\hskip 3pt X_4=\{u\},$$
primitive collections $X_0\cup X_1$, $X_1\cup X_2,\dots,X_4\cup
X_0$ and primitive relations:
$$v_1+\dots+v_{n-r}+y-cz_2-\dots-cz_{r}-(b+1)t=0,$$
$$y+z_1+\dots+z_{r}-u=0,$$
$$z_1+\dots+z_{r}+t=0,$$
$$t+u-y=0,$$
$$u+v_1+\dots+v_{n-r}-c_2z_2-\dots-c_{r}z_{r}-bt=0,$$
where $b$ and $c$ are positive integers.

Then the ordered collection of line bundles
$$L_{s,q}=\o(-sD_t-sD_y+(q-bs-(n-r))D_v)$$
for $s=0,\dots,r$ and $q=0,\dots,n-r$ and
$$L'_{s,q}=\o(-sD_t-(s-1)D_y+(q-bs-(n-r))D_v)$$
for $s=0,\dots,r-1$ and $q=0,\dots,n-r-1$ where the order
is defined by $L_{s,q}<L'_{s,q}<L_{s,q+1},\;L_{s+1,q_1}<L_{s,q_2}$
is a full, strongly exceptional collection of line bundles. \ktw
\dow From Subsection \ref{ssAcyc} we already know that this is a
strongly exceptional collection. We have just checked the sufficient condition for fullness in
Proposition \ref{13}.\kdow

\section{Bondal's construction not containing a full, strongly exceptional collection}\label{no}

\subsection{Example}

Let us consider the case when: $$X_0=\{v_1\},\quad
X_1=\{y_1,\dots,y_k\},\quad X_2=\{z_1\},$$$$
X_3=\{t_1,\dots,t_k\},\quad X_4=\{u_1,\dots,u_k\}$$ then we can take
$$v_1,y_2,\dots,y_k,t_1,\dots,t_k,u_2,\dots,u_k$$ to be a basis of
the lattice $N=\mathbbm Z^{3k-1}$. Other vectors are like in
\ref{rownania} with all coefficients $b_i$ and $c_i$ equal to zero. We have linear dependencies of divisors:
$$D_{v_1}=D_{u_1}+D_{y_1},\;\;\;D_{t_i}=D_{z_1}+D_{y_1},\;\;\;D_{y_i}=D_{y_1},\;\;\;D_{u_i}=D_{u_1}$$

Let $B$ be the image of the real torus in the Picard group as described in the Subsection \ref{ssBondThom}. One can easily see that:

$$B=\{\o([\sum_{i=1}^k-\alpha_t^i]D_{z_1}+[\sum_{i=2}^k-\alpha_u^i-\alpha_v^1]D_{u_1}+[-\alpha_v^1+\sum_{i=2}^k-\alpha_y^i+\sum_{i=1}^k\alpha_t^i]D_{y_1}):$$$$0\leq\alpha_v^i,\alpha_y^i,\alpha_t^i,\alpha_u^i<1\}.$$

So $B$ is contained in the set:
$$S:=\{\o(-aD_{z_1}-bD_{u_1}+(a-c)D_{y_1}):a,b,c\in\{0,\dots,k\}\}=$$
$$=\{\o(-a(D_{z_1}-D_{y_1})-bD_{u_1}-cD_{y_1}):a,b,c\in\{0,\dots,k\}\}.$$
From Corollary \ref{acyclic2} we know that line bundle is acyclic if and
only if it is not isomorphic to any of the following line bundles
$$\o(\alpha_1D_{v_1}+\alpha_2D_{y_1}+\alpha_3D_{z_1}+\alpha_4D_{t_1}+\alpha_5D_{u_1})=$$
$$=\o((\alpha_3+\alpha_4)(D_{z_1}-D_{y_1})+(\alpha_1+\alpha_2+\alpha_3)D_{y_1}+(\alpha_1+\alpha_5)D_{u_1}),$$
where exactly $2,3$ or $5$ consecutive $\alpha$ are negative and
if $\alpha_2<0$ then $\alpha_2\leq -k$, if $\alpha_4<0$ then
$\alpha_4\leq -k$ and if $\alpha_5<0$ then $\alpha_5\leq -k$. Let us
observe that line bundles form the set
$$R=\{\o(a(D_{z_1}-D_{y_1})+bD_{y_1}+cD_{u_1}):(a,b,c)\in[\frac{k}{2},k]\times[-k,-\frac{k}{2}-1]\times[0,k]\}$$
are not acyclic. Indeed fixing $\alpha_1=-k,\alpha_3=\frac{k}{2}$
and taking $\alpha_4,\alpha_5$ nonnegative and $\alpha_2$ negative
we can achieve all of them. Let us define the set of pairs
$$P:=\{-(\frac{k}{2}+\frac{a}{2})(D_{z_1}-D_{y_1})-(\frac{k}{2}+\frac{b}{2})D_{y_1}-(\frac{k}{2}+\frac{c}{2})D_{u_1},-(\frac{k}{2}-\frac{a}{2})(D_{z_1}-D_{y_1})-$$
$$-(\frac{k}{2}-\frac{b}{2})D_{y_1}-(\frac{k}{2}-\frac{c}{2})D_{u_1}):(a,b,c)\in[\frac{k}{2},k]\times[-k,-\frac{k}{2}-1]\times[0,k]\}.$$
It is easy to see that elements of these pairs are distinct and
they belong to $S$. Difference in each pair is an element of $R$
so it is not acyclic line bundle. Hence to have a strongly
exceptional collection $C$ in $S$ we have to exclude at least one
element from each pair. To have integer coefficients of divisors
in $P$ we should take $a\equiv b\equiv c\equiv k$ (mod $2$), so we
have to throw out at least $\frac{k^3}{32}$ elements among
$(k+1)^3$ elements in $S$. Full, strongly exceptional collection
has to have $l$ elements, where $l$ is the rank of the Grothendick
group $K^{0}(X)$ (for toric varieties it is isomorphic to $\mathbb
Z^{l}$, where $l$ is the number of maximal cones). In our case
there are at least $k^3$ maximal cones, since each time we throw
out one element from $X_2,X_4$ and $X_5$ we get different maximal
cone (exact number is $k^3+2k^2+2k$). So we have proven the
following:

\tw If $(k+1)^3-\frac{1}{32}k^3<k^3+2k^2+2k$, what is when $k>32$,
then there is no full, strongly exceptional collection among line
bundles that come from Bondal's construction. \ktw
\uwa
Notice that the considered variety is Fano, so is expected to have a full, strongly exceptional collection.
\kuwa

\subsection{Our case}

Let us consider the case from Subsection \ref{ssour}, but with all coefficients $c_i$ equal to $c\leq b$. Let $B$ be the image of the real torus in the Picard group as described in the Subsection \ref{ssBondThom}. One can see that:

$$B=\{\o([\sum_{i=1}^{r}-\alpha_z^i]D_{t}+[\sum_{i=1}^{n-r}-\alpha_v^i+c\sum_{i=2}^{r}\alpha_z^i -(b+1)\sum_{i=1}^{r}\alpha_z^i)]D_{y}+$$ $$+[\sum_{i=1}^{n-r}-\alpha_v^i+c\sum_{i=2}^{r}\alpha_z^i -b\sum_{i=1}^{r}\alpha_z^i)]D_{u}): 0\leq\alpha_v^i,\alpha_z^i<1\}.$$

So $B$ is contained in the set:
$$S:=\{\o(-sD_{t}-sD_{y}+qD_{v}),\o(-sD_{t}-(s-1)D_{y}+qD_{v}):s\in\{0,\dots,r\},$$
$$q\in\{-(n-r)-c-(b-c)s),\dots,(b-c)(-s+1)\}\}$$

Our collection defined in Subsection \ref{ssour}, or its torsion,
is contained in the set $S$ unless $cr\leq b$. It can be also
shown that if this inequality fails then there is no full strongly
exceptional collection among line bundles that come from Bondal's
construction.

\text{Micha\l\hskip 4pt Laso\'{n}}\linebreak
\text{Mathematical Institute of the Polish Academy of Sciences}\linebreak
\text{\'{S}w. Tomasza 30, 31-027 Krak\'{o}w, Poland}\linebreak
\text{Theoretical Computer Science Department}\linebreak
\text{Faculty of Mathematics and Computer Science}\linebreak
\text{Jagiellonian University, 30-348 Krak\'{o}w, Poland}\linebreak
\texttt{mlason@op.pl}
\vskip 5pt
\text{Mateusz Micha\l ek}\linebreak
\text{Mathematical Institute of the Polish Academy of Sciences}\linebreak
\text{\'{S}w. Tomasza 30, 31-027 Krak\'{o}w, Poland}\linebreak
\text{Institut Fourier, Universite Joseph Fourier}\linebreak
\text{100 rue des Maths, BP 74, 38402 St Martin d'H\'eres, France}\linebreak
\texttt{wajcha2@poczta.onet.pl}
\end{document}